\documentclass[10pt]{article}
\usepackage{latexsym,amsfonts,amssymb,amsmath,amsthm,graphicx}
        \newtheorem{theorem}{Theorem}
        \newtheorem{definition}{Definition}
        \newtheorem{proposition}{Proposition}[section]
        \newtheorem{corollary}{Corollary}[section]
        \newtheorem{lemma}{Lemma}

        \newtheorem{example}{Example}

\begin{document}

\title{The Algebra P_n is Koszul}
\author{David Nacin}

\abstract{The algebras $Q_n$ describe the relationship between the
roots and coefficients of a non-commutative polynomial. I.Gelfand,
S.Gelfand, and V. Retakh have defined quotients of these algebras
corresponding to graphs. In this work we find the Hilbert series
of the class of algebras corresponding to the $n$-vertex path,
$P_n$. We also show this algebra is Koszul.

We do this by first looking at class of quadratic algebras we call
Partially Generator Commuting.  We then find a sufficient
condition for a PGC-Algebra to be Koszul and use this to show a
similar class of PGC algebras, which we call ch$P_n$, is Koszul.
Then we show it is possible to extend what we did to the algebras
$P_n$ although they are not PGC.

Finally we examine the Hilbert Series of the algebras $P_n$}

\section{Koszul Algebras}

There are a number of equivalent definitions of Koszul algebras
including this lattice definition from Ufnarovskij \cite{U}.

\begin{definition} \label{AVR}

A quadratic algebra $A = \{V,R\}$ (where $V$ is the span of the
generators and $R$ the span of the generating relations in
$V\otimes V$) is Koszul if the collection of $n-1$ subspaces
$\{V^{\otimes i-1} \otimes R \otimes V^{\otimes n-i-1}\}_i$
generates a distributive lattice in $V^{\otimes n}$ for any $n$.

\end{definition}

The characterization of Koszulity we will need arises from this
definition and is based on the diamond lemma.  Suppose that $A$ is
a quadratic algebra with relations $R$ in $V \otimes V$ and a
monomial ordering exists so that every overlap ambiguity of degree
three resolves.  Then $A$ is a PBW-algebra (see chapter one of
\cite{PP}). The following result is due to S. Priddy and found in
\cite{Priddy}.

\begin{theorem} \label{PBW}

Any quadratic PBW-algebra is Koszul.

\end{theorem}

We will also need the following theorem from \cite{U}.

\begin{theorem} \label{dualkoszul}

A quadratic algebra $A$ is Koszul iff its dual algebra $A^*$ is
Koszul.  In the situation where they are both Koszul the Hilbert
series of $A$ is given by $\frac{1}{h(-x)}$ where $h(x)$ is the
Hilbert series of $A$.

\end{theorem}

\section{$Q_n$ and Some Related Algebras}

Let $P(x) = x^n - a_{n-1}x^{n-1} + a_{n-2}x^{n-2} - \cdots +
(-1)^na_0$ be a polynomial over a division algebra.  I. Gelfand
and V. Retakh \cite{vieta} studied relationships between the
coefficients $a_i$ and a generic set $\{x_1, \cdots, x_n\}$ of
solutions of $P(x)=0$. For any ordering $(i_1, \cdots, i_n)$ of
$\{1, \cdots, n\}$ one can construct \textit{pseudoroots} $y_k$,
$k = 1, \cdots n$, (certain rational functions in $x_{i_1},
\cdots, x_{i_n}$) that give a decomposition $P(t) = (t - y_n)
\cdots (t - y_2)(t - y_1)$ where $t$ is a central variable.

In \cite{GRW} I. Gelfand, V. Retakh, and R. Wilson introduced the
algebra $Q_n$ of all pseudo-roots of a generic noncommutative
polynomial, determined a basis for this algebra and studied its
structure. These algebras are quadratic and perhaps most easily
presented by generators $r(A)$ for all nonempty $A \subset [n] =
\{1, ... , n \}$ and relations
$$r(A)(r(A \setminus \{i\}) - r(A \setminus \{ j \})) + (r(A
\setminus \{ i \}) - r(A \setminus \{ j \}))r(A \setminus \{i,j
\}) - r(A \setminus \{i \})^2 + r(A \setminus \{j \})^2$$ for all
$i,j \in A \subset \{1,2, ... ,n\}$ where $r(\emptyset) = 0$. For
example, $Q_1$ is the free algebra with one generator (isomorphic
to $k[x]$) and $Q_2$ is the algebra with generators $r(1), r(2),
r(1,2)$ and the one relation $r(1,2)(r(1) - r(2)) = r(1)^2 -
r(2)^2$.  Though this definition is fairly straightforward, it is
the next presentation that will be more useful in our
construction.

The algebras $Q_n$ have a presentation given by generators $u(A),
\emptyset \neq A \subset [n]$ and relations

$$\sum_{C,D \subset A} [u(C \cup i), u(D \cup j)] = (\sum_{E \subset
A} u(E \cup i \cup j)) \sum_{F \subset A}(u(F \cup i) - u(F \cup
j))$$ for all $A \subset [n], i,j \in [n]\setminus A, i \neq j$.

\begin{definition} A complex with $n$ nodes is a family
$\cal F$ of nonempty subsets $A \subset [n]$ satisfying $A \in
{\cal F}, B \subset A \Rightarrow B \in {\cal F}$.  The dimension
of $\cal F$ is defined as $dim {\cal{F}} = max_{A \in \cal F}(|A|
- 1)$
\end{definition}

\begin{definition} Let ${\cal F}$ be a complex with $n$ nodes. Define $Q_n({\cal
F})$ to be the quotient algebra of $Q_n$ by the ideal generated by
the elements $u(A)$ for all $A \notin \cal F$.
\end{definition}

Notice that for any complex $\cal F$, $Q_n({\cal F})$ has a
presentation as a quadratic algebra.

\begin{example} If ${\cal F} = {\cal P}([n])- \emptyset$ then $Q_n({\cal F}) \cong
Q_n$ (Here ${\cal P}([n])$ denotes the power set, or collection of
all subsets of $[n]$.)
\end{example}

\begin{example} If ${\cal F}=\{A \subset {\cal F} | |A| =
1\}$ then $Q_n({\cal F})$ is isomorphic to the algebra of
commutative polynomials in $n$ variables.
\end{example}

If ${\cal F}' \subset {\cal F}$ is a subcomplex then $Q_n({\cal
F}')$ is naturally isomorphic to a quotient algebra of $Q_n({\cal
F})$.

Let $n_1 < n_2$ and let $\cal F$ be a complex with $n_1$ nodes.
Then as $[n_1] \subset [n_2]$, $\cal F$ may be viewed as a complex
with $n_2$ nodes.  We denote this complex by ${\cal F}'$.  Then
$Q_{n_1}({\cal F}) \cong Q_{n_2}({\cal F}')$ since every generator
$u(A)$ of $Q_{n_2}$ with $A \nsubseteq [n_1]$ is outside ${\cal
F}'$.  Consequently every algebra $Q_n({\cal F})$ occurs, up to
isomorphism, for a complex ${\cal F}$ containing every ${i}$, $1
\leq i \leq n$.

Consider the case where $\cal F$ is a complex of dimension one
with $n$ nodes. We can then also look at $\cal F$  as a graph on
$n$ nodes. To do this define $V$, our set of vertices, to be the
set of elements of $\cal F$ with cardinality one.  Our set of
edges, $E$ is the set of elements of $\cal F$ with cardinality
two. We adopt the convention of considering a graph to have no
loops or multiple edges.  Then there is actually a one to one
correspondence between graphs on $n$ vertices and complexes with
$n$ nodes and dimension one.

\begin{theorem} \label{graphalg} \cite{GGR} Let $\cal F$ be a complex with $n$ nodes and dimension one.
Then the algebra $Q_n({\cal F})$ is generated by the elements
$u(i)$ for $i \in [n]$ and $u(i,j)$ for $\{i,j\} \in E$ with the
following relations (assume $u(i,j)=0$ if $\{i,j\} \notin E$):

\noindent{$(i) [u(i),u(j)] = u(i,j)(u(i)-u(j)) \ i \neq j,\ i,j
\in [n]$}

\noindent{$(ii)
[u(i,k),u(j,k)]+[u(i,k),u(j)]+[u(i),u(j,k)]=u(i,j)(u(i,k)-u(j,k))$
for distinct $i,j,k \in [n]$}

\noindent{$(iii) [u(i,j),u(k,l)] = 0$ for distinct $i,j,k,l \in
[n]$}

\end{theorem}

If ${\cal F}$ is the complex of dimension one corresponding to a
graph $G$, we will write $Q_n(G) = Q_n({\cal F})$.  We refer to
the elements $u(i)$, $i \in [n]$, as nodes and $u(i,j)$, $\{i,j\}
\in E$, as edges.  We also often denote $u(i,j)$ as $u(ij)$. Using
this terminology the following proposition is immediate from (i)
and (iii).

\begin{proposition} Nodes in $G$ commute if they are not connected by an edge.
Non-adjacent edges in $G$ commute.
\end{proposition}

It is harder to find a way to simplify relation $(ii)$.  To gain
some insight first fix distinct $i,j,k$ in $[n]$.  Let $V= $
span$\{u(i), u(j), u(k), u(i,j), u(j,k), u(i,k)\}$ and let
$v_{i,j,k} =
[u(i,k),u(j,k)]+[u(i,k),u(j)]+[u(i),u(j,k)]-u(i,j)(u(i,k)-u(j,k))$.
Consider the natural action of $S_3$ (the permutation group on
three letters) on $T(V)$ defined by setting $\sigma u(i) =
u(\sigma (i))$ and $\sigma u(i,j) = u(\sigma(i), \sigma(j))$ and
extending linearly.

\begin{proposition} The orbit of $v_{i,j,k}$ under the action of $S_3$ spans a
space of dimension two in $T(V)$.
\end{proposition}

\begin{proof}  Let $\mu$ be the transposition given by $i \rightarrow j
\rightarrow i$ and $\tau$ be given by $i \rightarrow k \rightarrow
i$.  Since these permutations generate $S_3$ it will be enough to
show the action of $\tau$ and $\mu$ sends the space
span$\{v_{i,j,k}, v_{k,j,i} \}$ back to itself. Since this space
is clearly fixed by $\tau$ we need only worry about $\mu$. A short
computation shows that $\mu$ sends $v_{i,j,k}$ to $- v_{i,j,k}$
and $v_{k,j,i}$ to $v_{k,j,i} - v_{i,j,k}$.
\end{proof}

We say $G$ is triangle free if $\{i,j\}, \{j,k\} \in E
\Longrightarrow \{i,k\} \notin E$ for any $\{i,j\} \neq \{j,k\}$.
In this case relation $(ii)$ simplifies further.

\begin{proposition} \label{trifree} If $G$ is a triangle free graph then $(ii)$ is equivalent to

\noindent{$(ii')$  $u(i)$ commutes with $u(j,k)$ whenever $\{j,k\}
\in E, \{i,j\}, \{i,k\} \notin E$} and

\noindent{$(ii'') [u(i),u(jk)]+u(ij)u(jk) = [u(k), u(i,j)] +
u(j,k) u(i,j) = 0$ whenever $\{i,j\}, \{j,k\} \in E, \{i,k\}
\notin E $}

\noindent{for any distinct $i,j,k \in [n]$.}

\end{proposition}

\begin{proof} We know from our last proposition that we can replace
relation $(ii)$ with $v_{i,j,k}$ and $v_{k,j,i}$.  In the
situation where $\{i,j\}, \{i,k\} \notin E$, $v_{i,j,k}$ becomes
$[u(i), u(j,k)]$ and $v_{k,j,i}$ becomes zero since $u(i,j)$ and
$u(i,k)$ are zero. This gives us the relation $(ii')$.

In situation $(ii'')$ $u(i,k) = 0$ so $v_{i,j,k} =
[u(i),u(jk)]+u(ij)u(jk)$ and $v_{k,j,i} = [u(k), u(ij)] + u(jk)
u(ij)$ so we are done.
\end{proof}

\section{The Algebra $P_n$}

Now let us specialize to one particular triangle free graph, the
$n$ vertex path $P_n$ given in hypergraph notation by the complex
$$\{\{1\},\{2\}, \cdots, \{n\}, \{1,2\}, \{2,3\}, \cdots,
\{n-1,n\}\}.$$  To make our notation simpler we will identify
vertices and edges of the graph with the corresponding generators
for the algebra. Thus we refer to the elements of our algebra
$P_n$(which is really $Q_n(P_n)$) by $v_i$ for $u(i)$ and $e_{ij}$
for $u(i,j)$. Applying everything we have shown about the
relations $(i)$, $(ii)$ and $(iii)$ we get the following
proposition.

\begin{proposition} \label{Pn} The algebra $P_n$ generated by the complex (graph) $P_n$ is
presented by generators $v_1, v_2, \cdots, v_n, e_{12}, e_{23},
\cdots, e_{n-1,n}$ and relations \label{Pn}

$[v_i,v_j]=0$ for $j>i+1, \ i,j \in [n]$

$[v_i,v_{i+1}]+e_{i,i+1}(v_{i+1}-v_i)=0$ for $i \in [n-1]$

$[e_{i,i+1}, e_{j,j+1}]=0$ for $j>i+1, \ i,j \in [n-1]$

$[v_i,e_{j,j+1}]=0$ if $j>i+1$ or $j<i-2, \ j \in [n-1], i \in
[n]$

$[v_i,e_{i+1,i+2}]+e_{i,i+1}e_{i+1,i+2}=0, \ i \in [n-2]$

$[v_{i+2},e_{i,i+1}]+e_{i+1,i+2}e_{i,i+1}=0, \ i \in [n-2]$
\end{proposition}

Our goal in this chapter will be to show this algebra has the
Koszul property and find a way to compute its Hilbert series.

\section{$ch(P_n)$}

Let $V$ be the span of the generators of $P_n$.  We start by
defining an increasing filtration of $T(V)$.

\begin{proposition} \label{Pfilt}

Set $G^{(0)} = \textbf{F}1$. Then defining $G^{(i)}$ = span
$\{u(A_1)u(A_2) \cdots u(A_k) | \sum_{l=1}^k (3-|A_l|) \leq i\}$
for $i \geq 1$ defines a filtration of $T(V)$.

\end{proposition}

\begin{proof}

To show when $i \leq j$ that $G^{(i)} \subset G^{(j)}$ notice that
if $\sum_{l=1}^k (3-|A_l|) \leq i$ then $\sum_{l=1}^k (3-|A_l|)
\leq j$.  To show $\cup G^{(i)} = T(V)$ notice each monomial
$u(A_1)u(A_2) \cdots u(A_k)$ in $T(V)$ is contained in $G^{(i)}$
for $i = \sum_{l=1}^k (3-|A_l|)$.  Finally we must show that
$G^{(i)}G^{(j)} \subset G^{(i+j)}$.  To do this we will show that
the product of a monomial in $G^{(i)}$ and a monomial in $G^{(j)}$
is contained in $G^{(i+j)}$ then by extending linearly we will
know this is true for sums of monomials.

Suppose that $u(A_1)u(A_2) \cdots u(A_k) \in G^{(i)}$ and
$u(B_1)u(B_2) \cdots u(B_j) \in G^{(j)}$.  Then the product of
these two monomials is $u(A_1)u(A_2) \cdots u(A_k)u(A_{k+1})
\cdots u(A_{k+j})$ where $A_{k+i} = B_i$. Then $\sum_{l=1}^{k+j}
(3-|A_l|) = (\sum_{l=1}^{k} (3-|A_l|)) + (\sum_{l=1}^{j}
(3-|B_l|)) \leq i + j$.

\end{proof}

\begin{example}

Under this filtration $G^{(1)}$ is the span of all our $u(i,j)$.
This is because $(3-|\{i,j\}|) = 1$ so our monomial can only be of
length one.  Because $u(i)$ has $(3-|\{i\}|) = 2$ we get that
$G^{(2)}$ contains monomials of the form $u(i)$, $u(i,j)$, or
$u(i,j)u(k,l)$.

\end{example}

This induces a filtration of our algebra $P_n$ and hence we can
consider the associated graded algebra $gr(P_n)$.  If we chop off
the non-commutator terms in the relations given in proposition
\ref{Pn} those relations will all hold true in $gr(P_n)$. However,
there is no reason at this time to think that these are enough to
present $gr(P_n)$.  This does not stop us from considering the
algebra given by these chopped relations.  We call it $ch(P_n)$.

\begin{definition} The algebra $ch(P_n)$ is presented by generators $v_1, v_2,
\cdots, v_n, e_{12}, e_{23}, \cdots, e_{n-1,n}$ and relations

$[v_i,v_j]=0$ for $j>i+1, \ i,j \in [n]$

$[v_i,v_{i+1}]=0$ for $i \in [n-1]$

$[e_{i,i+1}, e_{j,j+1}]=0$ for $j>i+1, \ i,j \in [n-1]$

$[v_i,e_{j,j+1}]=0$ if $j>i+1$ or $j<i-2, \ j \in [n-1], i \in
[n]$

$[v_i,e_{i+1,i+2}]=0, \ i \in [n-2]$

$[v_{i+2},e_{i,i+1}]=0, \ i \in [n-2]$
\end{definition}

Our intention is to use Bergman's diamond lemma \cite{B} to find a
basis for $ch(P_n)$ and later $gr(P_n)$. Right now we can only see
that $gr(P_n)$ is a quotient of $ch(P_n)$.  We will soon show that
the two are indeed equal. In order to show equality in such a
situation it is enough to show that these two algebras have the
same Hilbert series. Before applying the diamond lemma to
$ch(P_n)$ we shall first prove some more general results about
algebras whose relations are all commutators of generators.  Then
we can see how these results will apply here.

\section{PGC-algebras}

\begin{definition} Suppose an algebra $A$ has a presentation by generators
$\{a_1,a_2, \cdots ,a_n\}$ and some relation set $R$ where each
relation in $R$ is of the form $[a_i,a_j]=0$ for some $i,j \in
[n]$. We then call $A$ a pre-generator-commuting algebra (or a
PGC-algebra for short).
\end{definition}


Notice that if $A$ is a PGC-algebra on n generators then for some
ideal $I$ in $A$, $A/I$ = $S_n$ the free commutative algebra on
$n$ generators.  However not every such algebra is PGC as the next
example shows.

\begin{example}

The algebra $A$ presented by generators $a,b$ and relation
$[ab,a]=0$ is not PGC.

\end{example}

\begin{example}

The algebra $A$ presented by generators $a,b,c$ and relation
$[a,b]=0$ is PGC.

\end{example}

\begin{example}

Both the commutative algebra and free algebra on $n$ generators
are PGC.

\end{example}

The structure of a PGC-algebra is based entirely on which pairs of
generators commute.  We can then describe such an algebra using a
graph with a node to represent each generator and edges between
two nodes if the generators commute.

\begin{definition}  Let $A$ be a PGC-algebra on generators $\{a_1,a_2, \cdots
,a_n\}$. The commuting graph $G_c(A)$ of $A$ is the graph with
nodes labeled $\{a_1,a_2, \cdots ,a_n\}$ and edges given by the
rule $e_{a_i,a_j}$ is an edge if and only if $[a_i,a_j]=0$.  We
often also consider the complement of this graph,
$\overline{G_c(A)}$. We call this the non-commuting graph of $A$
since an edge exists between two generators only when they do not
commute.
\end{definition}

\begin{example}  A simple example of a PGC-algebra is the algebra presented
by generators $\{a,b,c\}$ and relations $ab-ba=bc-cb=0$.  The
graph $G_c(A)$ is shown here:

\begin{center}

\includegraphics[height = .1\textwidth]{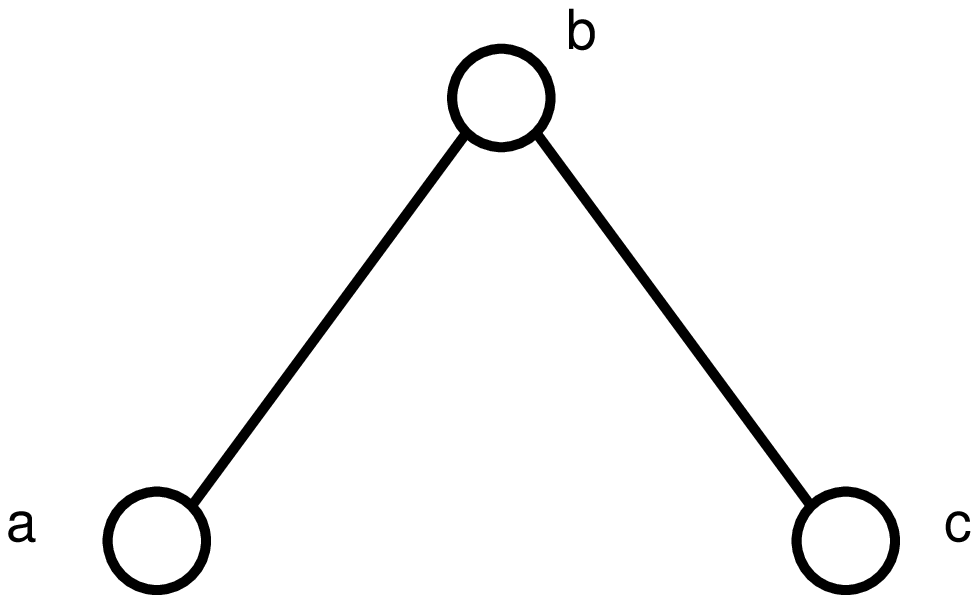}
\end{center}

\bigskip

The non-commuting graph $\overline{G_c(A)}$ is:

\begin{center}

\includegraphics[height = .1\textwidth]{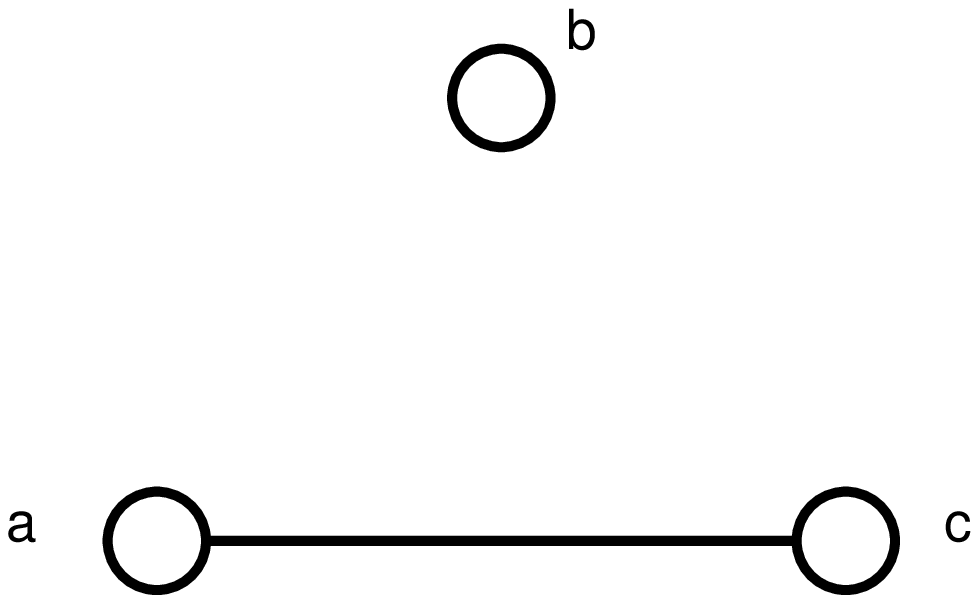}
\end{center}

\bigskip
\bigskip

The diamond lemma gives a method for determining a basis for such
an algebra. If, as in the example above, we choose a monomial
ordering given first by length and then lexicographically with
$c>b>a$ we get the following reductions:

$cb=bc$

$ba=ab$

This gives us one ambiguity, namely $cba$, that we must resolve.

$c(ba)=cab$ and $(cb)a=bac$ which gives us the new reduction
$cab=bac$ and one new ambiguity to resolve, $caba$.

$ca(ba)=caab$ and $(cab)a=baca$ which gives us the new reduction
$caab=baca$ and the new ambiguity $caaba$.

We can inductively show that by adjoining the reductions $caa
\cdots ab$ = $baca \cdots a$ we can resolve all ambiguities and we
end up with the following complicated list of bad words:
$cb,ba,cab,caab,caa \cdots ab$, $\cdots$.  A basis for the algebra
consists of the set of all monomials not containing one of these
strings.
\end{example}

\begin{example} Now let us look at the same algebra but apply the diamond
lemma with a different monomial ordering. First we order monomials
by length and then lexicographically with $b>a>c$. This gives us
the reductions:

$ba=ab$

$bc=cb$

This gives us no ambiguities and a basis for our algebra
consisting of all monomials not containing the strings ba or bc.
This is much simpler to use especially if we want to find the
Hilbert series of this algebra.
\end{example}

We are interested in these instances where all ambiguities of
degree three resolve (that is all ambiguities involving a monomial
of length three resolve) not only because is it easier to compute
the Hilbert series of such algebras.  Once we have shown that
there exists an ordering under the diamond lemma which causes
ambiguities of degree three to resolve, we can use theorem
\ref{PBW} to show our algebra is Koszul.  We will need one
definition and the following propositions, which are equivalent.

\begin{definition}  If $G$ is a graph with $n$ nodes then a vertex ordering of $G$ is
a surjective map from the vertices of $G$ onto $[n]$.
\end{definition}

\begin{proposition}  Let $A$ be a PGC-algebra with commuting graph ${G_c(A)}$.
Suppose there exists a vertex ordering of ${G_c(A)}$ so for any
three vertices $a, b$ and $c$ if $\{a,c\} \notin E,
\{a,b\},\{b,c\} \in E$ then neither $a<b<c$ nor $c<b<a$. Then
there exists a monomial ordering so that all ambiguities of degree
three are resolvable with the diamond lemma.
\end{proposition}

\begin{proposition} \label{nc} Let $A$ be a PGC-algebra with
non-commuting graph $\overline{G_c(A)}$. Suppose there exists a
vertex ordering of $\overline{G_c(A)}$ so for any three vertices
$a, b$ and $c$ if $\{a,c\} \in E, \{a,b\},\{b,c\} \notin E$ then
neither $a<b<c$ nor $c<b<a$.  Then there exists a monomial
ordering so that all ambiguities of degree three are resolvable
with the diamond lemma.
\end{proposition}

\begin{proof} Since the two statements are equivalent, we will prove only
the first.  Order monomials first by length and then by
lexicographically extending the vertex ordering.  It is enough to
show that given any three distinct vertices $a$,$b$ and $c$ that
all ambiguities involving those generators resolve.

First notice that if none of our generators $a$,$b$ and $c$
commute with each other, then there can be no ambiguity.  The same
holds if there is only one commuting pair.

If all three commute with each other, then we get one ambiguity
which is resolvable since everything commutes.

Finally, consider the case where $\{a,c\} \notin E,
\{a,b\},\{b,c\} \in E$.  The only ambiguities that could arise
would come from the monomials $abc$ or $cba$.  However, since $b$
is not in between $c$ and $a$ in the ordering it is not possible
for both $cb$ and $ba$ to be reductions (and similarly $ab$ and
$bc$). Hence there is no ambiguity to resolve here and we are
done.

\end{proof}

\begin{example}

The graph $P_4$ has three edges and four vertices:

\begin{center}

\includegraphics[height = .02\textwidth]{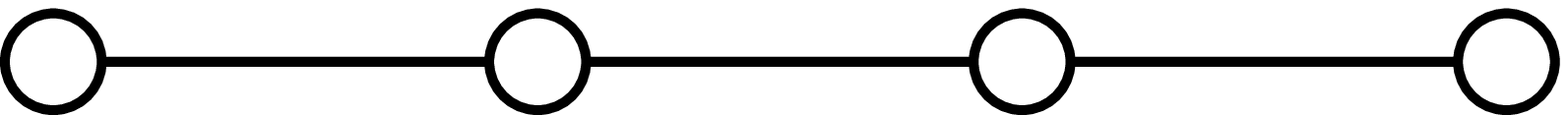}
\end{center}

The algebra $ch(P_4)$ is PGC hence we can take the non-commuting
graph. In the non-commuting graph of $ch(P_4)$ we have to
represent each generator with a vertex; this means one vertex for
each vertex in $P_4$ and one vertex for each edge in $P_4$.  What
we get looks like this (once we connect generators that do not
commute):

\begin{center}

\includegraphics[height = .1\textwidth]{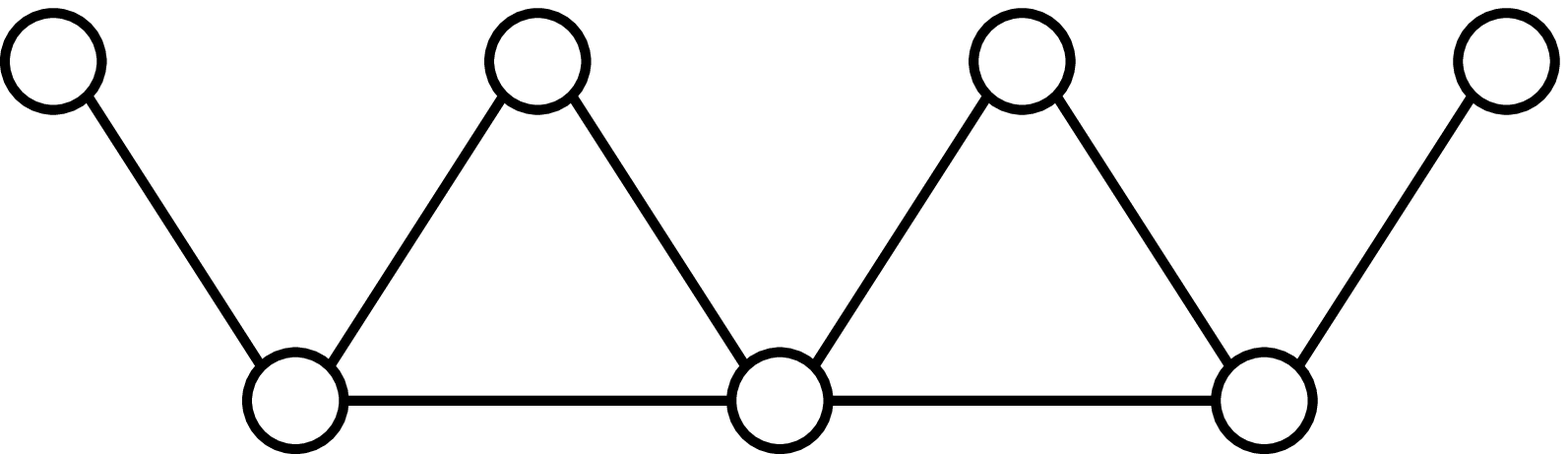}
\end{center}

We now need to find a labelling that satisfies the requirements of
proposition \ref{nc}.  The following labelling works:

\begin{center}

\includegraphics[height = .1\textwidth]{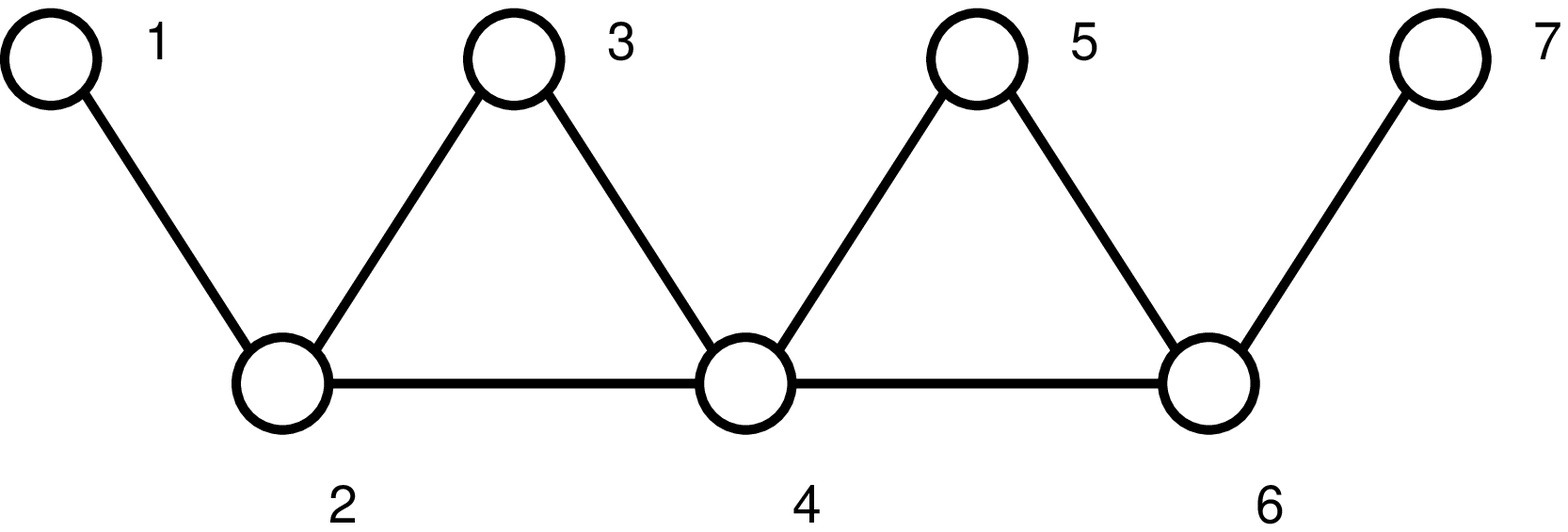}
\end{center}

With this vertex ordering we have shown the algebra's ambiguities
resolve.  This also tells us (by theorem \ref{PBW}) that $ch(P_4)$
is Koszul.

\end{example}

Of course, this technique does not only work for $ch(P_4)$

\begin{example}

From the graph

\begin{center}

\includegraphics[height = .02\textwidth]{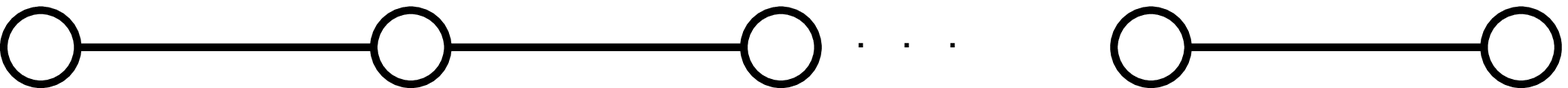}
\end{center}

$P_n$ we can form the non-commuting graph:

\begin{center}

\includegraphics[height = .1\textwidth]{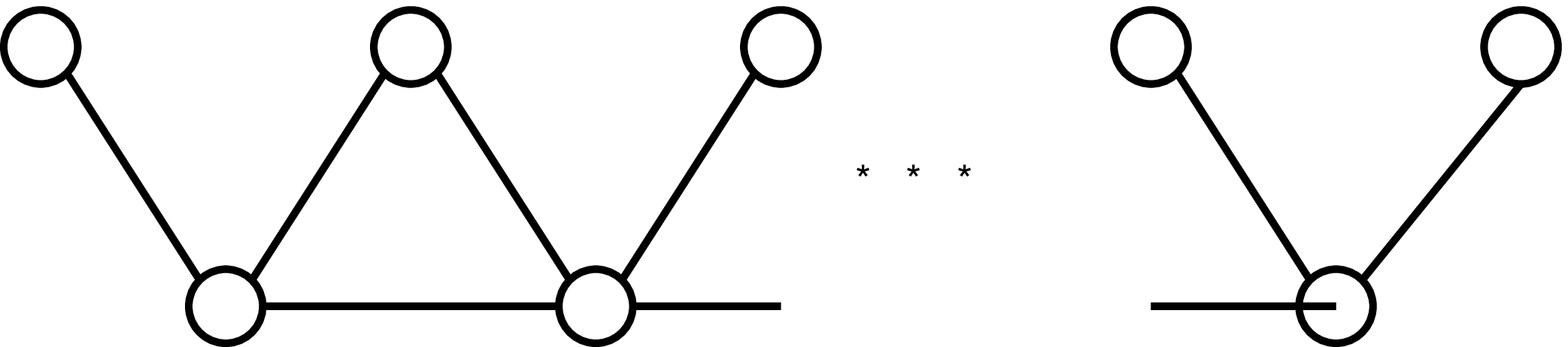}
\end{center}

Using proposition \ref{nc} and the ordering shown here

\begin{center}

\includegraphics[height = .1\textwidth]{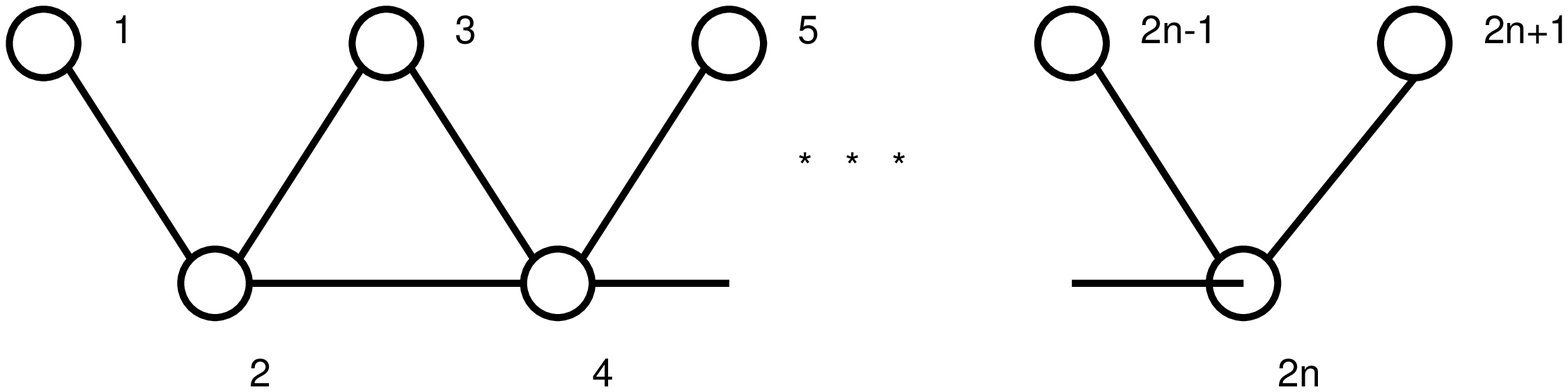}
\end{center}

we have shown the algebra's ambiguities resolve and by theorem
\ref{PBW}, $ch(P_n)$ is Koszul.

\end{example}

\section{Adaptation to $P_n$}

Now the propositions we have developed in this section only hold
for PGC-algebras.  This gave us an ordering of generators which
worked for $ch(P_n)$.  However, we do not know if this algebra has
anything at all to do with $P_n$.  What we have actually found is
an interesting ordering we can attempt on $P_n$.  We will use the
same ordering for $P_n$ and we will see that with this ordering
all ambiguities still resolve. This will show $P_n$ is Koszul,
describe a basis for $P_n$, and show that $P_n$ has the same
Hilbert series as $ch(P_n)$ implying that $ch(P_n) = gr(P_n)$.

\begin{lemma}

Consider the monomial ordering of $T(V)$ arising lexicographically
from the generator ordering with $v_n > e_{n-1,n} > v_{n-1} >
\cdots > e_{2,3} > v_2 > e_{1,2} > v_1$.  Then with the set of
reductions generated by the relations of $P_n$ under this
ordering, all ambiguities of degree three are resolvable.

\end{lemma}

\begin{proof}

Our reductions are

$v_kv_{k-1} \rightarrow v_{k-1}v_k - e_{k,k-1}v_{k-1} +
e_{k,k-1}v_k$ for $1 < k \leq n$

$v_kv_j \rightarrow v_jv_k$ for $1 \leq j < k \leq n$

$v_ke_{k-1,k-2} \rightarrow e_{k-1,k-2}v_k - e_{k,k-1}e_{k-1,k-2}$
for $2 < k \leq n$

$v_ke_{j,j-1} \rightarrow e_{j,j-1}v_k$  for $1 \leq j < k-1$ and
$k \leq n$

$e_{k,k-1}v_{k-2} \rightarrow
v_{k-2}e_{k,k-1}+e_{k-1,k-2}e_{k,k-1}$ for $2 < k \leq n$

$e_{k,k-1}v_j \rightarrow v_je_{k,k-1}$ for  $1 \leq j < k-2$ and
$k \leq n$

$e_{k,k-1}e_{j,j-1} \rightarrow e_{j,j-1}e_{k,k-1}$ for $1 < j <
k-1$ and $k \leq n$.

\noindent This gives us ambiguities of one these forms

$v_jv_kv_l$

$e_{j,j+1}v_kv_l$

$v_{j+1}e_{k,k+1}v_l$

$v_jv_ke_{l-1,l}$

$e_{j,j+1}v_ke_{l-1,l}$

\noindent for each of the four cases $j=k+1=l+2, j>k+1=l+2,
j=k+1>l+2, j>k+1>l+2$ and

$e_{j,j+1}e_{k,k+1}v_l$ for $j-1>k=l+1$ and $j-1>k>l+1$

$v_je_{k,k+1}e_{l,l+1}$ for $j-2=k>l+1$ and $j-2>k>l+1$

$e_{j,j+1}e_{k,k+1}e_{l,l+1}$ for $j-1>k>l+1$.

This gives us 29 cases which need to be checked for $P_n$.  The
first $n$ for which all 29 ambiguities actually appear is $n=7$
and by symmetry it is enough for $P_n$ to resolve the following
ambiguities in $P_7$: $v_5v_3v_1, v_4v_3v_1, v_4v_2v_1,
v_3v_2v_1,$ $v_4v_3e_{1,2},$ $v_5v_3e_{1,2},$ $v_5v_4e_{1,2},
v_6v_4e_{1,2},$ $v_4e_{3,2}v_1,$ $v_5e_{3,4}v_1,$ $v_5e_{2,3}v_1,$
$v_6e_{3,4}v_1,$ $e_{3,4}v_2v_1,$ $e_{4,5}v_2v_1,$
$e_{4,5}v_3v_1,$ $e_{5,6}v_3v_1,$ $e_{4,5}v_3e_{1,2},$
$e_{5,6}v_3e_{1,2},$ $e_{5,6}v_4e_{1,2},$ $e_{6,7}v_4e_{1,2}$
$e_{4,5}e_{2,3}v_1,$ $e_{5,6}e_{3,4}v_1,$ $e_{5,6}e_{2,3}v_1,$
$e_{6,7}e_{3,4}v_1,$ $v_5e_{3,4}e_{1,2},$ $v_6e_{3,4}e_{1,2},$
$v_6e_{4,5}e_{1,2},$ $v_7e_{4,5}e_{1,2},$ $e_{5,6}e_{3,4}e_{1,2}$.
Each of these is easily checked.

\begin{itemize}

\item $v_5v_3v_1$ is perhaps the easiest case.  If we begin by
applying a reduction to the last two terms we get $v_5v_3v_1 =
v_3v_5v_1 = v_3v_1v_5 = v_1v_3v_5 .$  These match, hence this
ambiguity is resolvable.

\item $v_4 ( v_3v_1 ) = v_4v_1v_3 = v_1v_4v_3 = v_1v_3v_4 -
v_1e_{3,4}v_3 + v_1e_{4,3}v_4$ and $( v_4v_3 ) v_1 = v_3v_4v_1 -
e_{3,4}v_3v_1 + e_{3,4}v_4v_1 = v_3v_1v_4 - e_{3,4}v_1v_3 +
e_{3,4}v_1v_4 = v_1v_3v_4 - v_1e_{3,4}v_3 + v_1e_{3,4}v_4 .$ This
shows that this ambiguity is resolvable.

\item $v_4 ( v_2v_1 ) = v_4v_1v_2 - v_4e_{2,1}e_1 + v_4e_{2,1}v_2
= v_1v_4v_2 - e_{2,1}v_4v_1 + e_{2,1}v_4v_2 = v_1v_2v_4 -
e_{2,1}v_1v_4 + e_{2,1}v_2v_4$ and $( v_4v_2 ) v_1 = v_2v_4v_1 =
v_2v_1v_4 = v_1v_2v_4 - e_{2,1}v_1v_4 + e_{2,1}v_2v_4 .$

\item This next case is longer but not really any more difficult.
$v_3 ( v_2v_1 ) = v_3v_1v_2 - v_3e_{2,1}v_1 + v_3e_{2,1}v_2 =
v_1v_3v_2 - e_{2,1}v_3v_1 + e_{3,2}e_{2,1}v_1 + e_{2,1}v_3v_2 -
e_{3,2}e_{2,1}v_2 = v_1v_2v_3 - v_1e_{3,2}v_2 + v_1e_{3,2}v_3 -
e_{2,1}v_1v_3 + e_{3,2}e_{2,1}v_1 + e_{2,1}v_2v_3 -
e_{2,1}e_{3,2}v_2 + e_{2,1}e_{3,2}v_3 - e_{3,2}e_{2,1}v_2 .$  If
we reduce the other way we get $( v_3v_2 ) v_1 = v_2v_3v_1 -
e_{3,2}v_2v_1 + e_{3,2}v_3v_1 = v_2v_1v_3 - e_{3,2}v_1v_2 +
e_{3,2}e_{2,1}v_1 - e_{3,2}e_{2,1}v_2 + e_{3,2}v_1v_3 = v_1v_2v_3
- e_{2,1}v_1v_3 + e_{2,1}v_2v_3 - v_1e_{3,2}v_2 -
e_{2,1}e_{3,2}v_2 + e_{3,2}e_{2,1}v_1 - e_{3,2}e_{2,1}v_2 +
v_1e_{3,2}v_3 + e_{2,1}e_{3,2}v_3.$  As these are equal, this
ambiguity is resolvable.

\item The next ambiguity is long but simple to resolve as well.
$v_4 ( v_3e_{2,1} ) = v_4e_{2,1}v_3 - v_4e_{3,2}e_{2,1} =
e_{2,1}v_4v_3 - e_{3,2}v_4e_{2,1} + e_{4,3}e_{3,2}e_{2,1} =
e_{2,1}v_3v_4 - e_{2,1}e_{4,3}v_3 + e_{2,1}e_{4,3}v_4 -
e_{3,2}e_{1,2}v_4 + e_{4,3}e_{3,2}e_{1,2} = e_{2,1}v_3v_4 -
e_{4,3}e_{2,1}v_3 + e_{4,3}e_{2,1}v_4 - e_{3,2}e_{1,2}v_4 +
e_{4,3}e_{3,2}e_{1,2} = e_{2,1}v_3v_4 - e_{2,1}e_{4,3}v_3 +
e_{2,1}e_{4,3}v_4 - e_{3,2}e_{1,2}v_4 + e_{4,3}e_{3,2}e_{1,2}$ and
$( v_4v_3 ) e_{2,1} = v_3v_4e_{2,1} - e_{4,3}v_3e_{2,1} +
e_{4,3}v_4e_{2,1} = v_3e_{2,1}v_4 - e_{4,3}e_{2,1}v_3 +
e_{4,3}e_{3,2}e_{2,1} + e_{4,3}e_{2,1}v_4 = e_{2,1}v_3v_4 -
e_{3,2}e_{2,1}v_4 - e_{2,1}e_{4,3}v_3 + e_{4,3}e_{3,2}e_{2,1} +
e_{2,1}e_{4,3}v_4$.

\item $v_5 ( v_3e_{2,1} ) = v_5e_{2,1}v_3 - v_5e_{3,2}e_{2,1} =
e_{2,1}v_5v_3 - e_{3,2}v_5e_{2,1} = e_{2,1}v_3v_5 -
e_{3,2}e_{2,1}v_5$ and $( v_5v_3 ) e_{2,1} = v_3v_5e_{2,1} =
v_3e_{2,1}v_5 = e_{2,1}v_3v_5 - e_{3,2}e_{2,1}v_5$.

\item $v_5 ( v_4e_{2,1} ) = v_5e_{2,1}v_4 = e_{2,1}v_5v_4 =
e_{2,1}v_4v_5 - e_{2,1}e_{5,4}v_{4} + e_{2,1}e_{5,4}v_{5}$ and $(
v_5v_4 ) e_{2,1} = v_4v_5e_{2,1} - e_{5,4}v_4e_{2,1} +
e_{5,4}v_5e_{2,1} = v_4e_{2,1}v_5 - e_{5,4}e_{2,1}v_4 +
e_{5,4}e_{2,1}v_5 = e_{2,1}v_4v_5 - e_{2,1}e_{5,4}v_4 +
e_{2,1}e_{5,4}v_5.$

\item $v_6 ( v_4e_{1,2} ) = v_6e_{1,2}v_4 = e_{1,2}v_6v_4 =
e_{1,2}v_4v_6$ and $( v_6v_4 ) e_{1,2} = v_4v_6e_{1,2} =
v_4e_{1,2}v_6 = e_{1,2}v_4v_6.$

\item $v_4 ( e_{2,3}v_1 ) = v_4v_1e_{2,3} + v_4e_{1,2}e_{2,3} =
v_1v_4e_{2,3} + e_{1,2}v_4e_{2,3} = v_1e_{2,3}v_4 -
v_1e_{3,4}e_{2,3} + e_{1,2}e_{2,3}v_4 - e_{1,2}e_{3,4}e_{2,3}$ and
$( v_4e_{2,3} ) v_1 = e_{2,3}v_4v_1 - e_{3,4}e_{2,3}v_1 =
e_{2,3}v_1v_4 - e_{3,4}v_1e_{2,3} - e_{3,4}e_{1,2}e_{2,3} =
v_1e_{2,3}v_4 + e_{1,2}e_{2,3}v_4 - v_1e_{3,4}e_{2,3} -
e_{1,2}e_{3,4}e_{2,3}.$

\item $v_5 ( e_{3,4}v_1 ) = v_5v_1e_{3,4} = v_1v_5e_{3,4} = v
_1e_{3,4}v_5 - v_1e_{4,5}e_{3,4}$ and $( v_5e_{3,4} ) v_1 =
e_{3,4}v_5v_1 - e_{4,5}e_{3,4}v_1 = e_{3,4}v_1v_5 -
e_{4,5}v_1e_{3,4} = v_1e_{3,4}v_5 - v_1e_{4,5}e_{3,4}.$

\item $v_5 ( e_{2,3}v_1 ) = v_5v_1e_{2,3} + v_5e_{1,2}e_{2,3} =
v_1v_5e_{2,3} + e_{1,2}v_5e_{2,3} = v_1e_{2,3}v_5 +
e_{1,2}e_{2,3}v_5$ and $( v_5e_{2,3} ) v_1 = e_{2,3}v_5v_1 =
e_{2,3}v_1v_5 = v_1e_{2,3}v_5 + e_{1,2}e_{2,3}v_5.$

\item $v_6 ( e_{3,4}v_1 ) = v_6v_1e_{3,4} = v_1v_6e_{3,4} =
v_1e_{3,4}v_6$ and $( v_6e_{3,4} ) v_1 = e_{3,4}v_6v_1 =
e_{3,4}v_1v_6 = v_1e_{3,4}v_6.$

\item $e_{3,4} ( v_2v_1 ) = e_{3,4}v_1v_2 - e_{3,4}e_{1,2}v_1 +
e_{3,4}e_{1,2}v_2 = v_1e_{3,4}v_2 - e_{1,2}e_{3,4}v_1 +
e_{1,2}e_{3,4}v_2 = v_1v_2e_{3,4} + v_1e_{2,3}e_{3,4} -
e_{1,2}v_1e_{3,4} + e_{1,2}v_2e_{3,4} + e_{1,2}e_{2,3}e_{3,4}$ and
$( e_{3,4}v_2 ) v_1 = v_2e_{3,4}v_1 + e_{2,3}e_{3,4}v_1 =
v_2v_1e_{3,4} + e_{2,3}v_1e_{3,4} = v_1v_2e_{3,4} -
e_{1,2}v_1e_{3,4} + e_{1,2}v_2e_{3,4} + v_1e_{2,3}e_{3,4} +
e_{1,2}e_{2,3}e_{3,4}.$

\item $e_{4,5} ( v_2v_1 ) = e_{4,5}v_1v_2 - e_{4,5}e_{1,2}v_1 +
e_{4,5}e_{1,2}v_2 = v_1e_{4,5}v_2 - e_{1,2}e_{4,5}v_1 +
e_{1,2}e_{4,5}v_2 = v_1v_2e_{4,5} - e_{1,2}v_1e_{4,5} +
e_{1,2}v_2e_{4,5}$ and $( e_{4,5}v_2 ) v_1 = v_2e_{4,5}v_1 =
v_2v_1e_{4,5} = v_1v_2e_{4,5} - e_{1,2}v_1e_{4,5} + e_
{1,2}v_2e_{4,5}.$

\item $e_{4,5} ( v_3v_1 ) = e_{4,5}v_1v_3 = v_1e_{4,5}v_3 =
v_1v_3e_{4,5} + v_1e_{3,4}e_{4,5}$ and $( e_{4,5}v_3 ) v_1 =
v_3e_{4,5}v_1 + e_{3,4}e_{4,5}v_1 = v_3v_1e_{4,5} +
e_{3,4}v_1e_{4,5} = v_1v_3e_{4,5} + v_1e_{3,4}v_{4,5}.$

\item $e_{5,6} ( v_3v_1 ) = e_{5,6}v_1v_3 = v_1e_{5,6}v_3 =
v_1v_3e_{5,6}$ and $( e_{5,6}v_3 ) v_1 = v_3e_{5,6}v_1 =
v_3v_1e_{5,6} = v_1v_3e_{5,6}.$

\item $e_{4,5} ( v_3e_{1,2} ) = e_{4,5}e_{1,2}v_3 -
e_{4,5}e_{2,3}e_{1,2} = e_{1,2}e_{4,5}v_3 - e_{3,2}e_{4,5}e_{1,2}
= e_{1,2}v_3e_{4,5} + e_{1,2}e_{3,4}e_{4,5} -
e_{2,3}e_{1,2}e_{4,5}$ and $( e_{4,5}v_3 ) e_{1,2} =
v_3e{4,5}e_{1,2} + e_{3,4}e_{4,5}e_{1,2} = v_3e_{1,2}e_{4,5} +
e_{3,4}e_{1,2}e_{4,5} = e_{1,2}v_3e_{4,5} - e_{2,3}e_{1,2}e_{4,5}
+ e_{1,2}e_{3,4}e_{4,5}.$

\item $e_{5,6} ( v_3e_{1,2} ) = e_{5,6}e_{1,2}v_3 -
e_{5,6}e_{1,2}e_{2,3} = e_{1,2}e_{5,6}v_3 - e_{1,2}e_{5,6}e_{2,3}
= e_{1,2}v_3e_{5,6} - e_{1,2}e_{2,3}e_{5,6}$ and $( e_{5,6}v_3 )
e_{1,2} = v_3e_{5,6}e_{1,2} = v_3e_{1,2}e_{5,6} =
e_{1,2}v_3e_{5,6} - e_{1,2}e_{2,3}e_{5,6}$

\item $e_{5,6} ( v_4e_{1,2} ) = e_{5,6}e_{1,2}v_4 =
e_{1,2}e_{5,6}v_4 = e_{1,2}v_4e_{5,6} + e_{1,2}e_{4,5}e_{5,6}$ and
$( e_{5,6}v_4 ) e_{1,2} = v_4e_{5,6}e_{1,2} +
e_{4,5}e_{5,6}e_{1,2} = v_4e_{1,2}e_{5,6} + e_{4,5}e_{1,2}e_{5,6}
= e_{1,2}v_4e_{5,6} + e_{1,2}e_{4,5}e_{5,6}.$

\item $e_{6,7} ( v_4e_{1,2} ) = e_{6,7}e_{1,2}v_4 =
e_{1,2}e_{6,7}v_4 = e_{1,2}v_4e_{6,7}$ and $( e_{6,7}v_4 ) e_{1,2}
= v_4e_{6,7}e_{1,2} =
 v_4e_{1,2}e_{6,7} = e_{1,2}v_4e_{6,7}.$

\item $e_{4,5} ( e_{2,3}v_1 ) = e_{4,5}v_1e_{2,3} +
e_{4,5}e_{1,2}e_{2,3} = v_1e_{4,5}e_{2,3} + e_{1,2}e_{4,5}e_{2,3}
= v_1e_{2,3}e_{4,5} + e_{1,2}e_{2,3}e_{4,5}$ and $( e_{4,5}e_{2,3}
) v_1 = e_{2,3}e_{4,5}v_1 = e_{2,3}v_1e_{4,5} = v_1e_{2,3}e_{4,5}
+ e_{1,2}e_{2,3}e_{4,5}.$

\item $e_{5,6} ( e_{3,4}v_1 ) = e_{5,6}v_1e_{3,4} =
v_1e_{5,6}e_{3,4} = v_1e_{3,4}e_{5,6}$ and $( e_{5,6}e_{3,4} ) v_1
= e_{3,4}e_{5,6}v_1 = e_{3,4}v_1e_{5,6} = v_1e_{3,4}e_{5,6}.$

\item $e_{5,6} ( e_{2,3}v_1 ) = e_{5,6}v_1e_{2,3} +
e_{5,6}e_{1,2}e_{2,3}
 =
 v_1e_{5,6}e_{2,3} + e_{1,2}e_{5,6}e_{2,3} = v_1e_{2,3}e_{5,6}
 + e_{1,2}e_{2,3}e_{5,6}$ and $( e_{5,6}e_{2,3} ) v_1  = e_{2,3}e_{5,6}v_1
 = e_{2,3}v_1e_{5,6} = v_1e_{2,3}e_{5,6} + e_{1,2}e_{2,3}e_{5,6}.$

\item $e_{6,7} ( e_{3,4}v_1 ) = e_{6,7}v_1e_{3,4} =
v_1e_{6,7}e_{3,4} = v_1e_{3,4}e_{6,7}$ and $( e_{6,7}e_{3,4} ) v_1
= e_{3,4}e_{6,7}v_1 = e_{3,4}v_1e_{6,7} = v_1e_{3,4}e_{6,7}.$

\item $v_5 ( e_{3,4}e_{1,2} ) = v_5e_{1,2}e_{3,4} =
e_{1,2}v_5e_{3,4} = e_{1,2}e_{3,4}v_5 - e_{1,2}e_{3,4}e_{4,5}$ and
$( v_5e_{3,4} ) e_{1,2} = e_{3,4}v_5e_{1,2} -
e_{3,4}e_{4,5}e_{1,2} = e_{3,4}e_{1,2}v_5 - e_{3,4}e_{1,2}e_{4,5}
= e_ {1,2}e_{3,4}v_5 - e_{1,2}e_{3,4}e_{4,5}.$

\item $v_6 ( e_{3,4}e_{1,2} ) = v_6e_{1,2}e_{3,4} =
e_{1,2}v_6e_{3,4} = e_{1,2}e_{3,4}v_6$ and $( v_6e_{3,4} ) e_{1,2}
=
 e_{3,4}v_6e_{1,2} = e_{3,4}e_{1,2}v_6 = e_{1,2}e_{3,4}v_6$

\item $v_6 ( e_{4,5}e_{1,2} ) = v_6e_{1,2}e_{4,5} =
e_{1,2}v_6e_{4,5} = e_{1,2}e_{4,5}v_6$ and $( v_6e_{4,5} ) e_{1,2}
=
 e_{4,5}v_6e_{1,2} = e_{4,5}e_{1,2}v_6 = e_{1,2}e_{4,5}v_6.$

\item $v_7 ( e_{4,5}e_{1,2} ) = v_7e_{1,2}e_{4,5} =
e_{1,2}v_7e_{4,5} = e_{1,2}e_{4,5}v_7$ and $( v_7e_{4,5} ) e_{1,2}
= e_{4,5}v_7e_{1,2} =
 e_{4,5}e_{1,2}v_7 = e_{1,2}e_{4,5}v_7.$

\item $e_{5,6} ( e_{3,4}e_{1,2} ) = e_{5,6}e_{1,2}e_{3,4} =
e_{1,2}e_{5,6}e_{3,4} = e_{1,2}e_{3,4}e_{5,6}$ and $(
e_{5,6}e_{3,4} ) e_{1,2} = e_{3,4}e_{5,6}e_{1,2} =
e_{3,4}e_{1,2}e_{5,6} = e_{1,2}e_{3,4}e_{5,6}.$

\end{itemize}

\end{proof}

\begin{corollary}

The algebra $P_n$ has a basis consisting of monomials in $T(V)$
not containing any of substrings (from $V \otimes V$) of the
following forms:

$v_kv_{k-1}$ for $1 < k \leq n$

$v_kv_j$ for $1 \leq j < k \leq n$

$v_ke_{k-1,k-2}$ for $2 < k \leq n$

$v_ke_{j,j-1}$ for $1 \leq j < k-1$ and $k \leq n$

$e_{k,k-1}v_{k-2}$ for $2 < k \leq n$

$e_{k,k-1}v_j$ for  $1 \leq j < k-2$ and $k \leq n$

$e_{k,k-1}e_{j,j-1}$ for $1 < j < k-1$ and $k \leq n$.

\end{corollary}

\begin{corollary}

$gr(P_n) \cong ch(P_n)$.

\end{corollary}

\begin{proof}

Since $P_n$ and $ch(P_n)$ have the same list of ``bad'' words they
have the same graded dimension.  $P_n$ has the same graded
dimension as $gr(P_n)$.  Since $gr(P_n)$ is a quotient of
$ch(P_n)$ with the same graded dimension, the two must be
isomorphic.

\end{proof}

\begin{theorem}

$P_n$ is Koszul.

\end{theorem}

\begin{proof}

All ambiguities of degree three resolve so by theorem \ref{PBW},
$P_n$ is Koszul.

\end{proof}

\section{Duals and Hilbert Series}

Now suppose $A$ is any PGC-algebra with generators $a_1,a_2,
\cdots , a_n$.  The relations will be $a_ia_j - a_ja_i$ for all
$\{a_i,a_j\} \in G_c(A)$.  Let $b_i$ be the linear functional
sending $b_i(a_i)=1, b_i(a_j)=0$ if $j \neq i$.  Then the dual
algebra $A^*$ is given by generators $b_1,b_2, \cdots , b_n$ and
relations

$b_ib_j+b_jb_i$ for all $i \neq j, \{a_i,a_j\} \in G_c(A)$,

$b_ib_j=0$ for all $i \neq j, \{a_i,a_j\} \notin G_c(A)$

${b_i}^2=0$ for all i

Notice that $A^*$ is a quotient of the exterior algebra $E(span
\{b_1, \cdots, b_n\})$ on $b_1,b_2, \cdots , b_n$.  A word in
$E(span \{b_1, \cdots, b_n\})$ is zero in $A^*$ if and only if it
contains two letters $b_i$ and $b_j$ so $\{a_i,a_j\} \notin
G_c(A)$ (and hence $b_ib_j$=0 in $A^*$).

We can use this to compute the Hilbert series of $P_n$.  Since we
have already shown that $ch(P_n)$ has the same Hilbert series as
$P_n$ we can work solely with $ch(P_n)$ and use the fact that it
is a PGC-algebra.  In fact, since $ch(P_n)$ is Koszul, our
strategy will be to compute the Hilbert series $H(x)$ of
$ch(P_n)^*$ and get $\frac{1}{H(-x)}$ for the Hilbert series of
$P_n$.

Remember that in $ch(P_n)$ our generators were $\{v_1, \cdots,
v_n,e_{1,2}, \cdots e_{n-1,n}\}$.  Call the span of these
generators $V$.  The relations were simply that all the $v_i$s
commute with each other and $e_{i,i+1}$ commutes with everything
except for $e_{i-1,i}, e_{i+1,i+2}, v_i, $ and $ v_{i+1}$. Let
$w_i (i \in [n])$ be the functional $V \longrightarrow \bf{C}$
defined by $w_i(v_i)=1, w_i(v_j)=0$ for $i \neq j$, and
$w_i(e_{j,j+1})=0$. Let $d_{i,i+1}, i \in [n-1],$ be the
functional defined by $d_{i,i+1}(e_{i,i+1})=1,
d_{i,i+1}(e_{j,j+1})=0$ if $j \neq i$, and $d_{i,i+1}(v_j)=0$.
Then, by the reasoning we developed for general PGC-algebras, we
get that $A^*$ is the algebra given by generators $\{w_1, \cdots,
w_n,d_{1,2}, \cdots d_{n-1,n}\}$ and relations

$d_{i,i+1}v_i=v_id_{i,i+1}=0$

$d_{i,i+1}v_{i+1}=v_{i+1}d_{i,i+1}=0$

$d_{i,i+1}d_{i+1,i+2}=d_{i+1,i+2}d_{i,i+1}=0$

\noindent together with relations stating that all generators
anti-commute. We will sometimes refer to the $w_i$s as nodes and
$d_{i,i+1}$s as edges, considering their origins.

Computing the Hilbert series of $ch(P_n)^*$ requires counting the
number of subsets $S$ of $\{w_1, \cdots, w_n, d_{1,2}, \cdots,
d_{n-1,n}\}$ so that for each $d_{i,i+1} \in S$ we know
$d_{i-1,i}, d_{i+1,i+2}, v_i, v_{i+1} \notin S$.

Suppose $S$ contains a total of $j$ $d$'s.  Notice that since no
two adjacent edges can be in $S$ the number of ways we can have
$j$ $d$'s is the number of matchings $M(n,j)$ of the graph $P_n$
of size $j$. Also, each edge rules out the possibility of exactly
two vertices. Hence if $|S|=i$ then we have ${n-2j \choose i-j}$
ways we can pick the vertices for $S$.  This gives us a total of
M(n,j)${n-2j \choose i-j}$ valid subsets containing $j$ $d$'s. The
total number of valid subsets $S$, with $|S|=i$ is then $\sum_{j
\geq 0} M(n,j){n-2j \choose i-j}$.  This tells us that if we write
the Hilbert series $H_n(x)$ of $ch(P_n)^*$ as $H_n(x) = H_n^0 +
H_n^1x + H_n^2x^2 + \cdots$ then we have $H_n^i = \sum_{j \geq 0}
M(n,j){n-2j \choose i-j}$. We set $H_n^i = 0$ in the cases where
$n<0$ or $i<0$.

$H_n^i$ is not the easiest thing to compute, so we will now find
some rules to make finding the coefficients of $H_n(x)$ easier.

\begin{proposition}

$H_n(x)$ is always a palindromic ($H_n^i = H_n^{n-i}$).

\end{proposition}

\begin{proof}

We want to show $H_n^i = H_n^{n-i}$ or that $\sum_{j=0}
M(n,j){n-2j \choose i-j} = \sum_{j=0} M(n,j){n-2j \choose n-k-j}$.
But since ${n-2j \choose i-j} = {n-2j \choose n-2j-k+j} = {n-2j
\choose n-j-k}$ we are done.

\end{proof}

If we write out a few $H_n(x)$ in a pyramid we notice that each
term is the sum of the three terms in the triangle above it.

$$x^2+3x+1$$
$$x^3+5x^2+5x+1$$
$$x^4+7x^3+13x^2+7x+1$$
$$x^5+9x^4+25x^3+25x^2+9x+1$$

We now establish this relation in general.  Since $H_n^k$ is the
coefficient of $x^k$ in $H_n(x)$, this is given by the following
proposition.

\begin{proposition}  For all $n \geq 2, k \geq 1$,
$H_n^k=H_{n-1}^k+H_{n-1}^{k-1}+H_{n-2}^{k-1}$.
\end{proposition}

\begin{proof}  We want

$$\sum_{j=0}M(n,j){n-2j \choose k-j} =$$

$$\sum_{j=0}M(n-1,j){n-1-2j \choose k-j} + \sum_{j=0}M(n-1,j){n-1-2j
\choose k-1-j} + \sum_{j=0}M(n-2,j){n-2-2j \choose k-1-j}$$

\noindent{}or simply

$$\sum_{j=0} ( M(n,j){n-2j \choose k-j}-M(n-1,j){n-1-2j \choose
k-j}-$$
$$M(n-1,j){n-1-2j \choose k-1-j} -M(n-2,j){n-2-2j \choose
k-1-j} )=0$$

\noindent{}Notice the expression we want to set to zero equals

$$\sum_{j=0}M(n,j){n-2j \choose k-j}-M(n-1,j)({n-1-2j \choose
k-j}+{n-1-2j \choose k-1-j}) -M(n-2,j){n-2-2j \choose k-1-j}$$

\noindent{}which by Pascal's identity is

$$\sum_{j=0}M(n,j){n-2j \choose k-j}-M(n-1,j){n-2j \choose k-j}-M(n-2,j){n-2-2j \choose k-1-j}$$

\noindent{}which is just

$$\sum_{j=0}(M(n,j)-M(n-1,j)){n-2j \choose k-j}-M(n-2,j){n-2-2j \choose k-1-j}$$

Notice also that a recurrence for the number of matchings of on
the line $P_n$ is given by $M(n,j)=M(n-1,j)+M(n-2,j-1)$.  Using
this on our first term gives us $$\sum_{j=0} M(n-2,j-1){n-2j
\choose k-j}-M(n-2,j){n-2-2j \choose k-1-j}$$ \noindent which
collapses to zero and we are done.
\end{proof}

By adding another variable we get the following result explaining
the coefficients of $H_n(x)$.

\begin{proposition} $\sum_{n \geq 0} t^n H_n(x) = \frac{1}{1 - t - xt - xt^2}$
\end{proposition}

\begin{proof} We want to show

$$\sum_{n \geq 0} (1 - t - xt - xt^2)t^nH_n(x) = 1$$

\noindent{}Since $H_n(x) = \sum_{k \geq 0} H_n^k x^k$ we have

$$\sum_{n \geq 0} (1 - t - xt - xt^2)t^nH_n(x)$$

$$= \sum_{n \geq 0} \sum_{k \geq 0} (1 - t - xt - xt^2)H_n^kt^nx^k$$

$$= \sum_{n,k \geq 0} (H_n^k - tH_n^k - xtH_n^k - xt^2H_n^k)t^nx^k$$

$$= \sum_{n,k \geq 0} H_n^kt^nx^k - \sum_{n,k \geq 0}
H_n^kt^{n+1}x^k - \sum_{n,k \geq 0} H_n^kt^{n+1}x^{k+1} -
\sum_{n,k \geq 0} H_n^kt^{n+2}x^{k+1}$$

$$= \sum_{n,k \geq 0} H_n^kt^nx^k - \sum_{n,k \geq 0} H_{n -
1}^kt^{n}x^k - \sum_{n,k \geq 0} H_{n-1}^{k-1}t^{n}x^{k} -
\sum_{n,k \geq 0} H_{n-2}^{k-1}t^{n}x^{k}$$

$$= \sum_{n,k \geq 0} (H_n^k - H_{n-1}^k - H_{n-1}^{k-1} -
H_{n-2}^{k-1}) x^kt^n$$

We wish to show that this sum is equal to one.  By our last
proposition we know that $(H_n^k - H_{n-1}^k - H_{n-1}^{k-1} -
H_{n-2}^{k-1}) = 0$ when $n \geq 2$ and $k \geq 1$.  We break the
sum into the three unknown cases of $k=0$; $n=0, k \geq 1$; and
$n=1, k \geq 1$ to get

$$\sum_{n \geq 0} (H_n^0 - H_{n-1}^0)t^n + \sum_{k \geq 1} H_0^k
x^k + \sum_{k \geq 1} (H_1^k - H_0^k - H_0^{k-1})tx^k$$

\noindent{}For this to equal one we need to show:

i) $H_0^0 = 1$,

ii)$H_0^k = 0$ for $k \geq 1$,

iii)$H_1^k - H_0^k - H_0^{k-1} = 0$ for $k \geq 1$ , and

iv)$H_n^0 - H_{n-1}^0 = 0$ when $n > 0$.

First $H_0^0 = \sum_{j \geq 0} M(0,j){0-2j \choose 0} = M(0,0){0
\choose 0} = 1$

For ii) notice $H_0^k = \sum_{j=0} M(0,j) {0 - 2j \choose k-j} =
 M(0,0){0  \choose k} = 0$ because $k \geq 1$.

By ii), statement iii) now reduces to showing $H_1^1 - H_0^0 = 0$
and $H_1^k = 0$ for $k \geq 2$.

Notice $H_1^k =\sum_{j=0} M(1,j) {1 - 2j \choose k-j} =
 M(1,0){1  \choose k}$.  Thus $H_1^k = 0$ for $k \geq 2$ and $H_1^1 = 1$.  Since $H_0^0$ is one by i), this case is done.

Finally, for iv) notice $H_n^0 = \sum_{j \geq 0} M(n,j){n-2j
\choose -j} = M(n,0){n \choose 0} = 1$ and thus $H_n^0 - H_{n-1}^0
= 0$ when $n > 0$.
\end{proof}

\bibliographystyle{plain}

\bibliography{citations}

\begin{thebibliography}{1}

\bibitem{B}
George~M. Bergman.
\newblock The diamond lemma for ring theory.
\newblock {\em Adv. in Math.}, 29(2):178--218, 1978.

\bibitem{GGR}
Israel Gelfand, Sergei Gelfand, and Vladimir Retakh.
\newblock Noncommutative algebras associated to complexes and graphs.
\newblock {\em Selecta Math. (N.S.)}, 7(4):525--531, 2001.

\bibitem{vieta}
Israel Gelfand and Vladimir Retakh.
\newblock Noncommutative {V}ieta theorem and symmetric functions.
\newblock In {\em The Gelfand Mathematical Seminars, 1993--1995}, Gelfand Math.
  Sem., pages 93--100. Birkh\"auser Boston, Boston, MA, 1996.

\bibitem{GRW}
Israel Gelfand, Vladimir Retakh, and Robert~Lee Wilson.
\newblock Quadratic linear algebras associated with factorizations of
  noncommutative polynomials and noncommutative differential polynomials.
\newblock {\em Selecta Math. (N.S.)}, 7(4):493--523, 2001.

\bibitem{PP}
Positselski~L. Polishchuk~A.
\newblock Quadratic algebras.
\newblock {\em Preprint}, 1996.

\bibitem{Priddy}
Stewart~B. Priddy.
\newblock Koszul resolutions.
\newblock {\em Trans. Amer. Math. Soc.}, 152:39--60, 1970.

\bibitem{U}
V.~A. Ufnarovskij.
\newblock Combinatorial and asymptotic methods in algebra [ {MR}1060321
  (92h:16024)].
\newblock In {\em Algebra, VI}, volume~57 of {\em Encyclopaedia Math. Sci.},
  pages 1--196. Springer, Berlin, 1995.

\end{thebibliography}

\end{document}